\documentclass[10pt]{amsart}
\usepackage{amsmath, amssymb}
\usepackage{graphicx}
\usepackage{subfigure}
 \usepackage{mathrsfs}
\newcommand{\no}[1]{#1}
\renewcommand{\no}[1]{}
\no{\usepackage{times}\usepackage[subscriptcorrection,
 slantedGreek, nofontinfo]{mtpro}
\renewcommand{\Delta}{\upDelta}}
\usepackage{color}

\usepackage[colorlinks,citecolor=red,pagebackref,hypertexnames=false,breaklinks]{hyperref}

\date{\today}
\newtheorem{lemma}{Lemma}[section]
\newtheorem{remark}{Remark}[section]

\newtheorem{corollary}{Corollary}[section]
\newtheorem{proposition}{Proposition}[section]
\newtheorem{theorem}{Theorem}[section]
\def\proof{{\bf Proof}  }
\def\eproof{\hspace{3cm}
    $\square $  \\}

\def\R{\mathbb{R}}
\def\C{\mathbb{C}}

\def\im{\textrm{Im}}
\def\re{\textrm{Re}}

\def\hat{\widehat}
\def\tilde{\widetilde}

\newcommand{\be}{\begin{equation}}
\newcommand{\ee}{\end{equation}}
\newcommand{\ba}{\begin{array}}
\newcommand{\ea}{\end{array}}

\newcommand{\bea}{\begin{eqnarray*}}
\newcommand{\eea}{\end{eqnarray*}}
\newcommand{\bean}{\begin{eqnarray}}
\newcommand{\eean}{\end{eqnarray}}

\def\cydot{\leavevmode\raise.4ex\hbox{.}}
\newcommand{\normmm}[1]{{\left\vert\kern-0.25ex\left\vert\kern-0.25ex\left\vert #1
    \right\vert\kern-0.25ex\right\vert\kern-0.25ex\right\vert}}

\newcommand\gobblepars{%
    \@ifnextchar\par%
        {\expandafter\gobblepars\@gobble}%
        {}}

\title[inverse medium]{Stability for the multifrequency
 inverse medium problem}

\author{ Gang Bao}

\address{Department of
Mathematics, Zhejiang University,  Hangzhou,
Zhejiang,  310027,    China}

\email{baog@zju.edu.cn}

\author{Faouzi Triki}

\address{Faouzi Triki,  Laboratoire Jean Kuntzmann,
UMR CNRS 5224, Universit\'e Grenoble-Alpes, 700 Avenue  Centrale,
38401 Saint-Martin- d'H\`eres, France}
\email{faouzi.triki@univ-grenoble-alpes.fr}

\begin{document}

\begin{abstract}
The  solution of a multi-frequency  1d inverse
medium  problem  consists of recovering
 the refractive index of a medium
from  measurements of the scattered
waves   for multiple frequencies. In this paper,
rigorous stability estimates are derived
when the frequency takes value in a bounded interval.
It is showed that the ill-posedness of the
inverse medium problem  decreases as the width of the
frequency interval becomes larger. More precisely, under
certain regularity assumptions on the refractive index,
the estimates indicate that the power in H\"older stability
is an increasing function of the  largest
value in the frequency  band. Finally, a Lipschitz
stability estimate is obtained for the observable part of the medium function
  defined through a truncated trace formula.

\end{abstract}

\maketitle
%\tableofcontents

\noindent {\footnotesize {\bf AMS subject classifications.} 35L05,
35R30, 74B05; Secondary 47A52, 65J20}

\noindent {\footnotesize {\bf Key words.}
Inverse medium problem, Helmholtz equation,  Stability estimates, trace formula, scattering resonances.}
%\newpage

%\tableofcontents

%\vfill\break

\section{Introduction}
This paper is concerned with  the stability  for determining
the refractive index of an one-dimensional (1d)  medium from boundary measurements. For a fixed
frequency, it is known that this inverse problem is severely ill-posed and suffers from the
 lack of uniqueness.
Several numerical results show that in the case of multiple frequencies, in contrast
 with the single  frequency case, the ill-posedness  decreases
dramatically when the frequency band increases and covers the
resonance region of the medium~(\cite{chen1},~\cite{BLiu}, ~\cite{BLi},~\cite{BHL} and
references therein). However, little is known about the stability for the inverse problem
or the convergence issues for the numerical methods. Our goal of the present paper is
to prove stability results for the multifrequency inverse medium scattering problem.
Such results would be essential for  a  rigorous justification
of the numerical observations.

Consider the 1d scalar Helmholtz equation
\bean \label{mainequation}
 \phi^{\prime \prime}(x,k)+k^2(1+q(x))\phi(x,k)  &=& 0,
\eean
where the real-valued $(1+q(x))$ is the refractive index of the medium.  For any real number
$k$, we look for a solution of  the form
\bea
\phi_\pm(x,k) &=& \psi_\pm(x,k)+e^{\pm ikx},
\eea
where the scattered wave $\psi_+,\psi_-$ corresponding to the left excitation
$e^{ikx}$, and the right excitation $e^{-ikx}$, respectively,  satisfy
the outgoing radiation conditions
$$
\begin{array}{llccc}
\psi^\prime(x,k)- ik\psi(x,k) &=& 0\qquad \qquad &\text{ for} \quad x\geq 1,\\
\psi^\prime(x,k) + ik\psi(x,k) &=& 0\qquad  \qquad &\textrm{ for}\quad x\leq 0.
\end{array}$$
The sum of the incident  wave and its corresponding  scattered wave, $\phi(x,k)$,
is called the total wave. Throughout, it is assumed  that  the medium function
$q(x)$ has the regularity
$ C^{m+1}_0\left([0,1]\right)$ with $m\geq 4$, and satisfies
\bean \label{refractivebounds}
 1+q(x) \geq n_0, \qquad \textrm{for }  x\in \R, \eean
with $n_0 \in (0, 1)$ a  fixed  constant.  The scattered wave $\psi(x,k)$ satisfies the Helmholtz equation
\bean \label{psi1}
\psi_\pm^{\prime \prime}(x,k)+k^2(1+q(x))\psi_\pm(x,k) &=&  -k^2q(x)e^{\pm ikx},
\label{psiequation}
\eean
for all $x\in (0,1)$.\\

Since $q(x) $ vanishes outside $(0,1)$, it is easy to see that for $k\in \mathbb R$, there exist
complex numbers $\mu_\pm$ known as the reflection coefficients, such that
\bean
\begin{array}{llccc}
\psi_+(x,k) &=& \mu_+(k)e^{-ikx} \qquad \textrm{for }\quad &x\leq 0,\\
\psi_-(x,k) &=& \mu_-(k)e^{ikx} \qquad \textrm{for }\quad &x\geq 1.
\end{array}
\label{muplusminus}
\eean
The  existence and uniqueness of the solutions $\psi_\pm \in C([0,1])$ are  well known
for any real $k$ \cite{BLT1}. \\

Therefore, the function $k\rightarrow \mu_\pm(k)$ are well defined on $\mathbb R$.
The outgoing radiation conditions imply
$$
\begin{array}{llccc}
\phi_+(x,k) &=& \phi_+(1,k)e^{ikx}  \qquad \textrm{for }\quad &x\geq 1,\\
\phi_-(x,k) &=& \phi_-(0,k)e^{-ikx} \qquad \textrm{for }\quad &x\leq 0.
\end{array}
$$
Furthermore, the constants $\phi_+(1,k)$ and $\phi_-(0,k)$ are nonzero. If
they are zero then Cauchy theorem implies that $\phi_\pm(x,k) =0 $ for all
$x\in \mathbb R$, which means that $\psi_\pm =-e^{\mp ikx}$ on the
whole space. This is in contradiction with the outgoing radiation
conditions. In fact,   $\phi_\pm(x,k)\not=0$ holds
for all $x\in \mathbb R$ and for all $k \in \mathbb C,$  satisfying
$ \im(k) \geq 0$ (Corollary 4.1~\cite{chen1}).  \\

The multifrequency inverse medium  problem may be stated as follows:\\

{\it Given one of  the reflection coefficients $\mu_+(k)$ and $\mu_-(k)$ for
 $ k\in (0, k_0) $, to reconstruct the
refractive  index $1+q(x)$ for $x\in [0,1]$.}\\

Define the impedance functions $p_\pm(x,k)$
associated with $\psi_\pm(x,k)$, respectively, by
\bean \label{impedancefunctions}
p_\pm(x,k)&=& \pm \frac{\phi_\pm^\prime(x,k)}{ik \phi_\pm(x,k)}.
\eean
It is shown in~\cite{chen1} that these functions are well defined
and verify  in addition the nonlinear Ricatti equation
\bean \label{impedanceequations}
p_\pm^\prime(x,k) \mp ik p_\pm^2(x,k) &=& \pm ik(1+q(x)),
\eean
subject to the boundary conditions
\bean
\ba{lllcccc} \label{impedanceboundaryconditions}
p_-(0,k) &=& 1;\; \; &   p_-(1,k) &=& d_-(k), \\
p_+(0,k)&=& d_+(k);\;\;&p_+(1,k)&=& 1,
\ea
\label{dplusminus}
\eean
for all $x\in (0,1)$, $k\in \mathbb R$, where
\bean\label{dfunctions}
d_\pm(k) &=& \frac{1-\mu_\pm(k)}{1+\mu_\pm(k)}.
\eean

The inverse problem may be restated
as: \\

{\it Given the data $d_-(k),\, k\in (0, k_0) $  or
$d_+(k),\, k\in (0, k_0), $
 to reconstruct the medium function
 $q(x)$ for $x\in [0,1]$. }\\

%%%%% discussion on Chen & Rokhlin's work
It is well known that in the case where the data is given for all
frequencies, this inverse problem
has a unique solution, and a number of algorithms have been proposed
for its numerical treatment~\cite{LS}. However, in applications,
the reflection coefficients  $\mu_\pm(k) $  are usually measured with  finite-accuracy
at a finite number of the frequencies $k$. Hence,
the well-posedness of the inverse problem  when
the measurements are taken over a finite interval is of critical importance. It is well known that the ill-posedness of the inverse
scattering problem decreases as the frequency increases \cite{ABFG}. However, at high frequencies,
the nonlinear equation becomes extremely oscillatory and possesses many more local minima.
A challenge for solving this problem is to develop a solution method that takes advantages of
the regularity of the problem for high frequencies without being undermined by local minima.
To overcome the difficulties, a recursive linearization method was proposed in \cite{chen1, chen2, chen3} for solving the inverse problem of the
two-dimensional Helmholtz equation. Based on the Riccati equations for the scattering
matrices, the method requires full aperture data and needs to solve a sensitivity matrix
equation at each iteration. The numerical results were very successful to address the ill-posedness computationally. However, there are two serious issues remain to be resolved.
Due to the high computational cost, it is numerically difficult to
extend the method to the three-dimensional problems. Recently, new and
more efficient recursive linearization methods have been developed  for solving the two-dimensional Helmholtz equation and the three dimensional
Maxwell equations for both full and limited aperture data by directly using the
differential equation formulation \cite{BLiu}, \cite{BLi}, \cite{BHL}, \cite{BLLT}. Theoretically, little is known about the stability for the inverse
problem with multiple frequency data. Our
main objective of this work is to establish stability estimates for the inverse problem with multiple frequency data.

%%%%%%%%%%%%%%%%%%%%%%%%%%%%%%%%%%%%%%%%%%%%
%%%%%%%%%%%%%%%%%%%%%%%%%%%%%%%%%%%%%%%%%%%%%
%\section{Main results} \label{mai}

We  state here our first main result associated to the inversion with  boundary measurements on
a band of frequencies.
  For $m \geq 4, \, M>0$,  and
$q_0 \in C^{m+1}_0([0,1])$ satisfying \eqref{refractivebounds},
 we further denote the set $\mathcal Q =
 \mathcal Q(n_0, m, M),$ by

 \bean \label{setQ}
\mathcal Q:= \{q \in C^{m+1}_0([0,1]): \; \|q-q_0\|_{C^{m+1}([0,1])} \leq M,\;  n_0\leq 1+q  \}.
\eean
We  next give our  first main stability estimate for the multifrequency inverse medium
problem.
In what follows
$c_{\mathcal Q}$ and $k_{\mathcal Q}$ denote generic  strictly positives constants depending only on
${\mathcal Q}$.
%%%%%%%%%%%%%%%%%%%%%%%%%%%

%%%%%%%%%%%%%%%%%%%%%%%%%%%%%%%
\begin{theorem}\label{mainhigh}
Assume that $q, \tilde q$ be two medium
functions in  $\mathcal Q$.  Let $d= d_\pm$ and $ \tilde d =\tilde d_\pm$
be  the boundary measurements associated respectively  to $q$
and $\tilde q$ as defined in~\eqref{dplusminus}, satisfying
$\|d -\tilde d \|_{L^\infty(0,+\infty)}  <1$. Let  $k^\star \in \mathbb R_+$ be the smallest
value
satisfying
 \bea
|d(k^\star) -\tilde d(k^\star)| = \|d -\tilde d \|_{L^\infty(0,+\infty)}.
\eea

 Then, there
exist constants $c_{\mathcal Q}>0$, and  $n_Q \in \mathbb N^*,$ such that
 the following estimate holds
 \bean \label{mainestimate}
 \left\| q- \tilde q \right\|_{L^\infty(\mathbb R)}
\leq c_{\mathcal Q}
 \|d -\tilde d \|_{L^\infty(0,k_0)}^{\frac{m}{m+1}w_0(k^\star, k_0)},
 \eean
 for all $k_0>0,$
where the function $w_0(k^\star, k_0)$ is continuous on $(\mathbb R_+^*)^2$, and  verifies
\bea
\frac{2}{\pi} \arctan(\frac{ (e^{k_0}-1)^{n_{\mathcal Q} }} {\sqrt{(e^{k^\star}-1)^{2n_{\mathcal Q}} - (e^{k_0}-1)^{2n_{\mathcal Q} } }})
\leq w_0(k^\star, k_0)\leq \frac{2}{\pi} \arctan\left(\inf\{\frac{k_0}{\sqrt{(k^\star)^2- k_0^2}}, \frac{ e^{k_0{n_{\mathcal Q}}}}{\sqrt{e^{2k^\star{n_{\mathcal Q}}}
- e^{2k_0{n_{\mathcal Q}}}} }\}\right),
\eea
for all $k_0 \in (0, k^\star]$.

\end{theorem}
%%%%%%%%%%%%%%%%%%%%%%%%%%%%
\begin{remark} \label{rem11}The H\"older exponent  $\frac{m}{m+1}w_0(k^\star, k_0)$
  in the estimate~\eqref{mainestimate} is an increasing function
of $k_0$.  It tends to zero when $k_0$ goes to zero which
shows as expected that the ill-posedness  of the inversion increases when
the band of frequency shrinks. On the other hand,
the  function  $w_0(k^\star, k_0)$ approaches to its upper bound $\frac{m}{m+1}$
when $k_0$ tends to $k^\star$, which is  the global H\"older stability estimate obtained
in Corollary~\ref{stabhigh}.  \\

 The value $k^\star$ represents the frequency  at which the noise is
the most important.  We observe that the  H\"older exponent $\frac{m}{m+1}w_0(k^\star, k_0)$ is
a decreasing function of $k^\star$, and tends to zero when $k^\star$ approaches  $+\infty$.\\

By considering  the stability estimate~\eqref{mainestimate}, we conclude  that the reconstruction
of the medium function
is accurate when the  frequency band  is large
enough and  contains the noise frequency ($k^\star \in (0, k_0]$), while it
 deteriorates when  the frequency band shrinks toward
zero.  These theoretical results confirm the numerical
observations and the physical expectations for the increasing
stability phenomena by taking multifrequency data.
\end{remark}

%%%%%%%%%%%%%%%%%%%%%%%%%%%%%%%%%%%%
\begin{theorem} \label{mainhigh3}
Assume that $q, \tilde q$ be two medium
functions in  $\mathcal Q$.  Let $d= d_\pm$ and $ \tilde d =\tilde d_\pm$
be  the boundary measurements associated respectively  to $q$
and $\tilde q$ as defined in~\eqref{dplusminus}, satisfying
$\varepsilon := \|d -\tilde d \|_{L^\infty(0,+\infty)}  <1$. \\

 Then, there
exist constants $c_{\mathcal Q}>0$, $k_{\mathcal Q}>0$ and  $n_Q \in \mathbb N^*,$ such that
 the following estimates hold
 \bean \label{mainestimate2}
 \left\| q- \tilde q \right\|_{L^\infty(\mathbb R)}
\leq c_{\mathcal Q} \varepsilon^{\frac{m}{m+1}}, \quad  \textrm{  if  }
k_0 \geq \frac{k_{\mathcal Q}}{ \varepsilon^{\frac{1}{m}}},\\ \label{mainestimate3}
\left\| q- \tilde q \right\|_{L^\infty(\mathbb R)}
\leq \frac{c_{\mathcal Q}}{ \left|\ln\left(\eta(k_0) |\ln(\varepsilon)|\right) \right|^{\frac{m^2}{m+1}}
}  \quad  \textrm{  if  }  k_0 < \frac{k_{\mathcal Q}}{ \varepsilon^{\frac{1}{m}}},
 \eean
where the function $\eta$ is given by
\bea
\eta(k_0) = \frac{(e^{k_0}-1)^{n_{\mathcal Q} } } {1+2\sqrt{1+ (e^{k_0}-1)^{2n_{\mathcal Q}}}}.
\eea
\end{theorem}
%%%%%%%%%%%%%%%%%%%%%%%%%%%%%%%%%%%%%
\begin{remark} The estimates \eqref{mainestimate2}  and \eqref{mainestimate3}
 show that the stability is H\"older when the  largest value in the frequency band $k_0$ is larger than a critical limit, and is of  logarithmic  type when $k_0$ becomes small. Hence for a limited band of frequencies one can improve the stability of the inverse problem by  increasing the largest frequency.  The critical limit  only depends on the noise in the measurement and the set  of medium functions $\mathcal Q$. When
 $k_0$ tends to zero the function $\eta(k_0)$ approaches zero, and right-hand side term blows up. This behavior  demonstrate 
 again  that the  inverse problem is severely  ill-posed  when $k_0$ is close to zero, and confirms the observations  made  in Remark \ref{rem11}.

 \end{remark}

Based on the
high frequency asymptotic expansions of the fields $\phi_\pm$,
~Chen and Rokhlin~\cite{chen1}  introduced the observable part
of the medium  $q(x)$ on the band of frequency $(0,k_0)$,
as the function $q_{k_0}(x)$  unique solution  to  the
truncated version of the trace formula~\eqref{partialtraceformula},
that is
\bean \label{truncated}
 p_{k_0, \pm}^\prime(x,k) \mp ik p_{k_0,\pm}^2(x,k) \mp ik(1+q(x)) &=&0,\\
\label{tracebis}
q^\prime(x) - \frac{2}{\pi}(1+q(x))
\int_{-k_0}^{k_0}(p_{k_0,+}(x,k)-p_{k_0,-}(x,k))dk &=&0,
\eean
for all $x\in (0,1)$, subject to the boundary conditions
\bean
\ba{lllccccccc} \label{boundarytruncated}
p_{k_0,+}(0,k) &=& d_+(k);\; \; &   p_{k_0, -}(0,k) &=& 1;\;\;  & q_{k_0}(0) = 0.
\ea
\eean
for all $k\in \mathbb C_+$. They also have derived error estimates
of the approximation of the medium function $q(x)$ by its observable
part $q_{k_0}(x)$ on the frequency band $(0, k_0)$ \cite{chen1}.\\

 Our third main result  is  to
 characterize $q_{k_0}(x)$ in terms of the
 frequency band $(0, k_0)$, and to show that  the recovery of  $q_{k_0}(x)$ is not sensitive
 to errors in the measurements if $k_0$ is large enough. \\
%%%%%%%%%%%%%%%%%%%%%%%%%%%%%
 %Based on the
%high frequency asymptotic expansions of the fields $\phi_\pm$,
%~Y. Chen and V. Rokhlin~\cite{chen1}  have defined the observable part
%of the medium  $q(x)$ on the band of frequency $(0,k_0)$,
%as the unique solution $q_{k_0}(x)$   to  the
%truncated version of a trace formula (which is provided in~\eqref{partialtraceformula}),
%that is
%\bean \label{rr1}
% p_{k_0, \pm}^\prime(x,k) \mp ik p_{k_0,\pm}^2(x,k) \mp ik(1+q(x)) &=&0,\\
%\label{rr2}
%q^\prime(x) - \frac{2}{\pi}(1+q(x))
%\int_{-k_0}^{k_0}(p_{k_0,+}(x,k)-p_{k_0,-}(x,k))dk &=&0,
%\eean
%for all $x\in (0,1)$, subject to the boundary conditions
%\bean\label{rr3}
%\ba{lllccccccc}
%p_{k_0,+}(0,k) &=& d_+(k);\; \; &   p_{k_0, -}(0,k) &=& 1;\;\;  & q_{k_0}(0) = 0.
%\ea
%\eean
%for all $k\in \mathbb C_+$.   They also derived error estimates
%of the approximation of the medium function $q(x)$ by its observable
%part $q_{k_0}(x)$ on the frequency band $(0, k_0)$ (see
%Lemma~\ref{qpapproximation} in Section \ref{obs}).
%We next  show that the observable part $q_{k_0}(x)$ is indeed stable under small
%perturbations in  boundary measurements.
%%%%%%%%%%%%%%%%%%%%%%%%%%%%%%%%%%
\begin{theorem}\label{estimateobservable}
Assume that $q, \tilde q$ be two medium
functions in  $\mathcal Q$.
 Let $d= d_\pm$ and $ \tilde d =\tilde d_\pm$
be  the boundary measurements associated respectively  to $q$
and $\tilde q$ as defined in~\eqref{dplusminus}.
 Let $q_{k_0}$ and $\tilde q_{k_{0}}$ be the observable
parts of respectively  $q$ and $\tilde q$ on $(0, k_0)$
solutions to the system~\eqref{truncated}-\eqref{tracebis}-\eqref{boundarytruncated}.
Then there exist constants
 $\rho_{\mathcal Q}>0$  and $k_{\mathcal Q}>0$
such that
\bea
\left\| q_{k_0}- \tilde q_{k_0} \right\|_{L^\infty(\mathbb R)}
 \leq \rho_{\mathcal Q}
\|d(k) -\tilde d(k)\|_{L^1(0, k_0)},
\eea
is satisfied  for all $k_0 \geq  k_{\mathcal Q}$.
\end{theorem}
%%%%%%%%%%%%%%%%%%%

%%%%%%%%%%%%%%%%%%%%%%%%%%%%%%%%%%%%%%%%%%%%%%%%

For higher dimension, to the best of our  knowledge, this inverse
 problem is still open. This is due to the difficulties in the analysis
of the scattering data  as a function of  the frequency, which are
related to  the strong nonlinearity  for high frequencies
and the existence of trapped rays.  From a physical point of view,
the situation is  better understood. According to
Uncertainty Principle  there exists a resolution limit to the
sharpness of details  of the medium that can be observed from
measurements in the far field region. This limit known as the
diffraction limit is about one half of the wavelength.  Consequently
the reconstruction of the medium can be then reduced by
increasing  the magnitude of the frequency~\cite{BT1}.
Mathematically, the inverse medium problem  with
 full  measurements at a fixed frequency is notoriously
ill-posed \cite{KV, SU}. In fact, Alessandrini  proved that the
stability estimates in 3d is of logarithmic type~\cite{Al}, and
Mandache  showed later the optimality of such estimates~\cite{Man}.
Recent studies have been conducted on the behavior of
the constant in the
logarithmic stability in terms of the fixed
frequency~\cite{ABFG, Isa, NUW}. Several other
results in inverse scattering  problems that are
related to the increasing stability phenomena
 by increasing the frequency were obtained
 in different settings~\cite{ACTV, ACZ, Isa, SN}.
 All of these results demonstrate the increasing stability
phenomena when the frequency becomes larger. For the case of
the inverse source problem for Helmholtz
equation and an homogeneous background it was shown
in~\cite{BLT1, BLT2, BLT3, CIL, IL1, IL2} that the ill-posedness of the inverse
problem decreases as the frequency increases. Convergence
results for iterative algorithms solving the multi-frequency
inverse medium problem  are obtained in~\cite{BT1,HQS}. Finally, we refer the reader to the
topical review on inverse scattering problems\cite{BLLT} with multifrequencies on other related topics.
%\textcolor{red}{
%In~\cite{BT1},
%the concept of the observable part of the medium for a given band of
%frequency is introduced.  Denote $q_{k_0}(x)$ the part of the
% medium function $q(x)$ that can  be recovered in a stable way from
% measurements for $k\in (0, k_0)$. Our second main result of this paper summarized
 %in Theorem~\ref{estimateobservable}, first  characterizes  $q_{k_0}(x)$ in terms of the
% frequency band $(0, k_0)$, and   shows that  the recovery of  $q_{k_0}(x)$
%  is insensitive to noise in the measurements if $k_0$ is large enough. }\\

The rest of the paper is structured as follows.  Auxiliary  results  related to the behavior
of the impedance functions as functions of the frequency are
provided in Section \ref{imp}.  The
 stability  estimate for the observable part of the medium is
 given in Section \ref{obs}.
 Finally, the proof of the main stability
estimates  for the multifrequency inverse
medium problem  is  provided in Sections \ref{proofmai} and \ref{proofmai2}.
%%%%%%%%%%%%%%%%%%%%%%
%%%%%%%%%%%%%%%%%%%%%%%%%%%%%%%%%%%%%%%%%%%%%%%%
\section{Properties of the impedance functions} \label{imp}
%%%%%%%%%%%%%%%%%%%%%%%%%%%%

A major difficulty in studying the multifrequency inverse medium problem is  the
fact that  the partial differential equation describing the scattering phenomena
involves a product of the frequency and the refractive index. In the 1d case,  Gel'fand-Levitan
techniques  can be employed  when the medium function is smooth
to convert the Helmholtz equation into a Schr\"odinger
equation. In the  obtained Schr\"odinger  equation, the frequency and the refractive index are separated,
which allows a   better understanding of the behavior of the solutions as functions of  the frequency.
This approach was used  to study the 1d  inverse spectral problem \cite{RS}.
It also led the  authors in~\cite{chen1} to derive  high-frequency asymptotic expansions
of the impedance functions. Here, we first present some of these useful asymptotic  results and our further analysis.
In addition, we also
 study   the meromorphic extensions of the  impedance functions
to the lower half complex plane.

For convenience, we complexify $k$. Denote $\mathbb C_\pm$ the upper half and lower half of the complex plane, that is
\bea
\mathbb C_+\;=\;\{k\in\mathbb C: \im(k) \geq 0\};\quad \mathbb C_-\;=\;\{k\in\mathbb C: \im(k) < 0\}.
\eea

It is easy to check from the uniqueness of the
equations~\eqref{impedanceequations} with the boundary
conditions~\eqref{impedanceboundaryconditions}, that
\bean \label{complexconjugatep}
 \overline{p_\pm(x,k)} &= &p_\pm(x,-\overline k),
\eean
for all $x\in \mathbb R$ and $k\in  \mathbb C_+$.\\

%%%%%%%%%%%%%%%%%%%%%%%%%%%%%%%
\subsection*{Low frequency behavior}
We next present  the behavior of the impedances functions when the frequency $k$ is close to
$0$.  In Lemma 4.1 and 4.2 of \cite{chen3},  the author derived the first
term in the asymptotic expansion  $p_\pm$  when $k$ approaches $0$. Here we provide
 explicit bounds in a given frequency neighborhood of $0$.

%%%%%%%%%%%%%%%%
\begin{proposition} \label{lowfrequency1}
 The following estimate
 \bea
  |d_\pm(k)| \leq 2,
 \eea
 holds for  all  $k \in \mathbb C$ satisfying $|k| \leq   1/M_1,$ with

 \bean \label{constantM1}
 M_1=   2(\|q_0\|_{L^\infty(0,1)}+M).
 \eean
\end{proposition}
%%%%%%%%%%%%%%
\proof
Since the proofs of the estimates for $d_+$ and $d_-$ are identical, we only
provide the proof for $d_+$. \\

Let
\[g_0(x,y)= \frac{e^{-ik|x-y|}}{-2ik},\]
be the Green function of the  one dimension Helmholtz equation with the same  radiation conditions
as $\psi_+$.  Multiplying the  equation \eqref{psi1} by $g_0(x,y)$ and integrating by parts yield the
following Lippmann-Schwinger integral equation
\bean \label{neumannseries}
(I_d-K_q)[\psi_+] = K_q[e^{-ik\cdot}],
\eean
where $I_d$ is the identity operator from $L^\infty(0,1)$ to itself, and $K_q$ is a linear  integral operator
on $L^\infty(0,1)$, defined by
\[K_q[\psi](x) = -k^2\int_0^1
g_0(x,y) q(y) \psi(y),\]
for all $\psi \in L^\infty(0,1)$.
Therefore for $2|k|(\|q_0\|_{L^\infty(0,1)}+M) \leq 1,$ the operator $K_q$ becomes a contraction, and we
deduce from the convergence of the Neumann series
\[
|\mu_+(k)| \leq \|\psi_+\|_{L^\infty(0,1)} \leq 1/3.
\]
Hence $|d_+(k)| \leq  2$ for $|k| \leq 1/M_1$,  which finishes the proof.

 \eproof

%%%%%%%%%%%%%%%%%%%%%%%
\begin{remark} \label{born} (Born approximation)
Using the Neumann series and after a forward calculation,  we obtain
\bea
\mu_+(k) &=& -\frac{k}{2i} \mathcal F(q)(-2k)+
\sum_{p=2}^\infty \left(\frac{ik}{2}\right)^p \int_{(0,1)^p} e^{ik \kappa_p(\xi)}
 Q_p(\xi)  d\xi,
\eea
for all $k\in (0,k_0)$, where $k_0<  1/M_1$,
$\kappa_p(\xi) = \xi_1+\sum_{l=1}^{p-1}|\xi_{l+1}-\xi_l| +\xi_p$
for all $\xi \in \mathbb R^p$, and $Q_p(\xi) =  \prod_{j=1}^pq(\xi_j)
$.  Since the first term in the low
frequency expansion is the Fourier transform  $\mathcal F(q)(2k),\,
k\in (-k_0,k_0)$, it seems natural to try to reconstruct
the medium function from this term by considering the
rest as a small perturbation ($O(k^2_0)$), and by using
 the same techniques as in~\cite{BLT1}. It turns out
that this approach fails to give any
approximation of the medium function. The Born approximation
error $O(k^2_0)$ is a higher order differential operator
that is exponentially amplified in the inversion
of the first term, and the final term does not vanish when
$k_0$ tends to zero.

\end{remark}

\subsection*{High frequency behavior} The following result was
obtained in~\cite{chen1}.

%%%%%%%%%%%%%%%%%%%%%%%%%%%%%%%
\begin{proposition} \label{thmppm}
 Assume that $q\in \mathcal Q$.
The impedances $p_\pm(x,k)$ are continuous functions of $(x,k)\in [0,1]\times \mathbb C_+$, and analytic
functions of $k \in \mathbb C_+$.
Moreover there exists a constant $c_{\mathcal Q}>0$ such that
the following estimates
\bean
\left \|p_\pm(x,k) - \sqrt{1+q(x)} \pm
 \frac{q^\prime(x)}{4i(1+q(x))} \frac{1}{k} \right \|_{L^\infty}
&\leq& \frac{c_{\mathcal Q}}{|k|^2}, \label{asympp} \\  \label{asymdiffIpm}
\left\| \overline{p_+(x,k)}-p_-(x,k)\right\|_{L^\infty}
&\leq& \frac{c_{\mathcal Q}}{|k|^m},
\eean
hold for all  $k\in \mathbb C_+^*$.
\end{proposition}
%%%%%%%%%%%%%%%%%%%%%%%%%%%%%%%%%
We remark that the estimate~\eqref{asympp}  provides
the two first terms in WKB expansions of the functions
$p_\pm$.  For large real $k$, the difference between
$\overline p_+$ and $p_-$ is extremely small, which decays
as $1/k^m$ where $m$ is the smoothness of the medium
$q(x)$.

%%%%%%%%%%%%%%%%%%%%%%%%%%
\subsection*{Meromorphic extension}  It is known that the
impedance functions $p_\pm(x, k) $ and
in particular the reflexion coefficients $\mu_\pm(k)$ are holomorphic
in $\mathbb C_+$, and have meromorphic extensions in  $\mathbb C_-$.
The poles of  $\mu_\pm$ are called the scattering resonances of the medium.
Here, we establish the existence
of a scattering resonances-free strip in the complex plane. The proof is based
on  a similar result for the 1d Schr\"odinger  equation derived in~\cite{Hit}.\\

From \eqref{mainequation} it follows that the poles can be characterized in the following
way: $k \in \C_-$ is a scattering pole if and only if there exists a nontrivial function $\phi$, such that

\bean \label{eigensystem1}
 \phi^{\prime \prime}(x,k)+k^2(1+q(x))\phi(x,k)  &=& 0, \quad x\in (0,1),
\eean
with
\bean\label{eigensystem2}
\phi^\prime(0,k) = -ik \phi(0,k), \;\; \phi^\prime(1) = ik \phi(1,k),
\eean
We now present a connection between the solution of the Helmholtz equation \eqref{eigensystem1}
 and the one of an equivalent  Schr\"odinger equation. This will  allow us to relate our scattering
resonances  to  the well studied poles of the resolvent of  the Schr\"odinger operator. This approach
has   been  also used  to derive the high frequency  asymptotic expansions  in Theorem~\ref{thmppm}.\\

Define  further the functions $n, x, N, r, \xi: \mathbb R \rightarrow \mathbb R$  by the following
expressions:
\bean
 n(x)= \sqrt{1+q(x)},\;
 t(x)= \int_0^x n(s) ds,\;
 N(t)= n(x(t))^{-1/4},\label{defNt}\\
 r(t) = \frac{N^{\prime\prime}(t)}{N(t)}-\frac{n^\prime(x)}{2(n(x))^2}=\frac{1}{4}n^{-4}(x)\left(q^{\prime\prime
 }(x)-nq^\prime(x)-\frac{5}{4}n^{-1}(q^\prime(x))^2\right).  \label{defr}
\eean
Then $\xi(t, k)$ defined by the Liouville transformation
\bea
 \xi(t, k):= N^{-1}(t)\phi(x(t),k),
 \eea
satisfies the Schr\"odinger equation:
\bean \label{eqpsi}
 \xi^{\prime \prime}(t, k) + (r(t)+k^2)\xi(t,k)  &=& 0,  \quad x\in (0,T),
\eean
with
\bean\label{eqpsi2}
\xi(0,k) = 1+\mu_+(k).
\eean
where $T = t(1)= \int_0^1 n(s) ds$ is the the travel time needed for the wave with speed $\frac{1}{n}$ to
propagate from one end to another.
 We remark that $r(t)$ has a compact support in $(0, T)$.
Consequently $k$ is a scattering resonance of~\eqref{eigensystem1}-\eqref{eigensystem2} iff  it is
a  resonance of the system~\eqref{eqpsi}-\eqref{eqpsi2}.  \\

The pole distribution of the resolvent for the  Schr\"odinger operator has been the subject of extensive investigations due to the continuous advance
of quantum mechanics. Many studies  have focused on the problem
of locating poles in the complex plane  for different  classes of potentials~\cite{Fr, Ha, Zw, De}.
For the one  dimensional  Schr\"odinger operator with super-exponentially decaying potentials,
 more precise results are possible. Particularly,   using the representation
of the  scattering matrix given by Melin~\cite{Me}, Hitrik~\cite{Hit} derived  an explicit  pole-free strip
for the Schr\"odinger operator  in the case of compactly supported potentials. The following result is
a direct consequence of  Hitrik's result and the observation that scattering resonances
 of the system~\eqref{eigensystem1}-\eqref{eigensystem2} are also the poles of  the Schr\"odinger operator~\eqref{eqpsi}-\eqref{eqpsi2}.
%%%%%%%%%%%%%%%%%%%%%%%%%%%%%%%%%%%%
\begin{proposition}\label{freestrip} Let $r:\; \mathbb R \rightarrow \mathbb R$ be defined by~\eqref{defr}, and
 $h(r):= \frac{1}{4T} e^{-2T\|r\|_{L^1(0,T)}}.$
Then the  strip
\bean
S_q = \left\{k\in \mathbb C;\, -h(r) \leq \im(k) \leq 0, \re(k)
\not= 0\right\}. \label{strip}
\eean
is free from scattering resonances of the system~\eqref{eigensystem1}-\eqref{eigensystem2}.
\end{proposition}

 %%%%%%%%%%%%%%%%
\begin{corollary} \label{remfreezone}
Let $c_{\mathcal Q, 1} = \max_{q\in \mathcal Q} \|n(x)\|_{L^\infty}$, and
$c_{\mathcal Q, 2} = \max_{q\in \mathcal Q}\|r(t)\|_{L^\infty}$. Then it follows from
Proposition~\ref{freestrip}  that the strip  of width $h_{\mathcal Q, 1}=  \frac{1}{4 c_{\mathcal Q, 1}}
e^{-2c_{\mathcal Q, 1}^2c_{\mathcal Q, 2}}$,
defined by
\bean
S_{\mathcal Q}^*= \left\{k\in \mathbb C;\, -h_{\mathcal Q,1} \leq \im(k) \leq 0, \re(k)
\not= 0\right\}. \label{strip2}
\eean
is free from scattering resonances of~\eqref{eigensystem1}-\eqref{eigensystem2} for all $q \in \mathcal Q$.
\end{corollary}
%%%%%%%%%%%%%%%%%%%%
We also deduce from  Proposition~\ref{freestrip} and Proposition~\ref{lowfrequency1}  that the  coefficients $d_\pm(k)$ have holomorphic
extensions in the strip $S_{\mathcal Q}$ defined by
 \bean \label{strip}
 S_{\mathcal Q}:= \left\{k\in \mathbb C;\, |\im(k)|  < h_{\mathcal Q} \right\},
\eean
where
\bea
h_{\mathcal Q} = \min \{h_{\mathcal Q,1}, \frac{1}{M_1}\}.
\eea

 We next obtain global bounds of these functions in the strip.
%%%%%%%%%%%%
\begin{proposition} \label{upperd}
There exist constants $k_{\mathcal Q}>0 $, $c_{\mathcal Q}>0$,
  $d_{\mathcal Q}>0$ that only depend on $Q$, such that
the following inequality  hold
\bean\label{dbounds2}
|d_\pm(k) -1| \leq \frac{c_{\mathcal Q}}{|\re(k)|^2},\quad  \forall k \in S_{\mathcal Q}, \; \re(k)  \geq k_{\mathcal Q},
\eean
\bean
\label{dbounds}
|d_\pm(k)| \leq d_{\mathcal Q},  \quad  \forall k \in S_{\mathcal Q}.
\eean

\end{proposition}
%%%%%%%%%%%%%
\proof
Since the proofs of the bounds for $\mu_+(k)$ and $\mu_-(k)$ are identical we
only provide the proof for the second  scattering coefficient. The
proof may be given by combining the general idea in the proof
of Lemma 4.12 in~\cite{chen2} and the meromorphic extension result above.

Applying the  Liouville  transformation to \eqref{mainequation}, we find
that  $\xi_-(t,k):= N^{-1}(t)\phi_-(x(t),k)$
satisfies the Schr\"odinger equation:
\bean \label{eqxi1}
 \xi^{\prime \prime}(t, k) + (r(t)+k^2)\xi(t,k)  &=& 0,  \quad t\in (0,T),
\eean
with
\bean\label{eqxi2}
\xi(t,k) = e^{ikt}, \quad t\leq 0.
\eean

The impedance function $p_-(x,k)$
is then given by
 \[p_-(x,k) = -n(x)\frac{\xi^\prime_-(t,k)}{ik\xi_-(t,k)}+ \frac{n^\prime(x)}{2ikn(x)}.\]

Introducing now the  auxiliary functions $\mathfrak m(t,k)= e^{ikt}\xi_-(t,k)$ and
$\mathfrak n(t,k)= -\frac{1}{ik}e^{ikt}\xi_-^\prime(t,k)$.\\

A forward calculation yields

\[
p_-(x,k) = \frac{\mathfrak m(t,k)}{\mathfrak n(t,k)}.
\]
We deduce from the system~\eqref{eqxi1}-\eqref{eqxi2}, that $\mathfrak m(t, k)$
satisfies
\bean \label{eqm1}
\mathfrak m^{\prime \prime}(t, k) -2ik \mathfrak m^\prime(t,k)  &=& -r(t)\mathfrak m(t,k),  \quad t\in (0,T),
\eean
with the initial conditions
\bean\label{eqm2}
\mathfrak m(0,k) = 1& \mathfrak m^\prime(0,k) = 0.
\eean
Multiplying \eqref{eqm1} by $e^{-2ikt}$ and integrating, we get
\bean \label{eqmm}
\mathfrak m^\prime(t,k) = -\int_0^t r(s)e^{2ik(t-s)}\mathfrak m(s, k) ds
\eean
Integrating the equation \eqref{eqmm}, we obtain
\bean \label{eqfredholm}
\mathfrak m &=& \frac{1}{2ik}\int_0^t  r(s)(1-e^{2ik(t-s)})\mathfrak m(s, k) ds+1,\\
 &=& \mathcal M_k[\mathfrak  m]+1,
\eean
where $\mathcal M_k: C(0,T)\rightarrow  C(0,T)$ is  a compact operator
defined by
\bean\label{mk}
\mathcal M_k[f](t) = \frac{1}{2ik}\int_0^t  r(s)(1-e^{2ik(t-s)})f(s) ds.
\eean

Since $q$ belongs to $\mathcal Q$ there exist constants $k_{\mathcal Q}>0$,
and $c_{\mathcal Q}>0$ such that
\bea
\|\mathcal M_k\|\leq \frac{c_{\mathcal Q}}{|\re(k)|}, \quad \forall k\in S_{\mathcal Q},\; |\re(k)|\geq k_{\mathcal Q}.
\eea

Then, the Fredholm equation~\eqref{eqfredholm} has a unique solution satisfying
\bea
|\mathfrak m(t,k)-1| \leq \frac{2c_{\mathcal Q}}{|\re(k)|}, \quad \forall t\in (0, T), \;
 \forall k\in S_{\mathcal Q},\; |\re(k)|\geq k_{\mathcal Q}.
\eea
It  can be approximated by the Neumann's series truncated at the second term
\bea
\mathfrak m(t,k) = 1+\frac{1}{2ik}\int_0^t r(s)ds+O(\frac{1}{|\re(k)|^2}),  \;
 \forall k\in S_{\mathcal Q},\; |\re(k)|\geq k_{\mathcal Q},
\eea
uniformly in $t \in (0, T)$.\\

Similarly, following the same approach, we have
\bea
\mathfrak  n(t,k) = 1+\frac{1}{2ik}\int_0^t r(s)ds+O(\frac{1}{|\re(k)|^2}),  \;
 \forall k\in S_{\mathcal Q},\; |\re(k)|\geq k_{\mathcal Q},
\eea
uniformly in $t \in (0, T)$.\\

Consequently
\bean \label{ppp}
\left|\frac{\mathfrak m(T,k)}{\mathfrak n(T,k)} -1\right|= |p_-(1,k)-1| = O(\frac{1}{|\re(k)|^2}),  \;
 \forall k\in S_{\mathcal Q},\; |\re(k)|\geq k_{\mathcal Q},
\eean
 Combining \eqref{ppp} with Proposition~\ref{lowfrequency1}, and
 the fact that $d_-(k)= p_-(1,k)$ is holomorphic in
$S_{\mathcal Q}$, we deduce the bound \eqref{dbounds} for $d_-$.\eproof

%%%%%%%%%%%%%%%
%\begin{remark}
%By calculating the two first terms in the Neumann series  of $\mathfrak m(t,k)$ and $\mathfrak n(t,k) $
 %we are able to prove that the asymptotic expansions in Proposition \ref{thmppm} are also valid
 %in the strip $S_{\mathcal Q}$.
%\end{remark}
%%%%%%%%%%%%%%%%%%%%%%%%%%%%%%%%%%%%%

%%%%%%%%%%%%%%%%%%%%%%%%%%%%%%%%%%%%%
\section{Observable part of the medium} \label{obs}
%\textcolor{red}{The concept of the observable part of the medium for a given band of
%frequency is introduced in~\cite{BT1}.
Recall from (\ref{truncated})-(\ref{boundarytruncated}) that the observable part of the
medium $q_{k_0}(x)$ for  $k\in (0, k_0)$. In this section using the truncated
 trace formula introduced in~\cite{chen1}, we
 characterize  $q_{k_0}(x)$ in terms of the
frequency band $(0, k_0)$, and study how its determination  is sensitive
to errors in the measurements.

The following trace formula is on the asymptotic behavior  in Proposition \ref{thmppm}.
%%%%%%%%%%%%%%%
\begin{lemma}(Trace formula, \cite{chen1})Let $q \in \mathcal Q$. Then the following  trace  formula holds
\bean \label{completetraceformula}
q^\prime(x) &=& \frac{2}{\pi}(1+q(x))
\int_{-\infty}^{\infty}(p_+(x,k)-p_-(x,k))dk.
\eean

More precisely, there exists a constant $c_{\mathcal Q}>0$ such that
the estimate
\bean \label{partialtraceformula}
\left \|q^\prime(x) - \frac{2}{\pi}(1+q(x))
\int_{-k_0}^{k_0}(p_+(x,k)-p_-(x,k))dk\right \| _{L^\infty(\mathbb R)} \leq
 \frac{c_{\mathcal Q}}{k_0^{m}},
\eean
holds for all $k_0\in \mathbb C^*$.

\end{lemma}
%%%%%%%%%%%%%%
The truncated version of the trace formula~\eqref{partialtraceformula} means that
the function
\bea
 \frac{2}{\pi}
\int_{-k_0}^{k_0}(p_+(x,k)-p_-(x,k))dk,
\eea
provides  a good approximation of  $\log(1+q(x))^\prime$ as long as
$k_0$ is large and the medium $q(x)$ is smooth.
%%%%%%%%%%%%%%%%%%%%%%%%%%%%%%%%%%%%%%%%%%%%%%%%%%%%%%%%%%%
\begin{lemma}  \label{qpapproximation} Let $q \in \mathcal Q$. Then, there exist constants $c_{\mathcal Q}>0$ and $k_{\mathcal Q}>0$
such that truncated trace formula
system~\eqref{truncated}-\eqref{tracebis}-\eqref{boundarytruncated}
 has a unique solution $q_{k_0}$. In addition the following estimates hold
\bea
\left\|p_{\pm}-p_{k_0, \pm} \right\|_{C\left( [0,1]\times[-k_0,k_0]\right)}, \;
\left\|q-q_{k_0} \right\|_{L^\infty(\mathbb R)} \leq \frac{c_{\mathcal Q}}{ k_0^{m}},
\eea
for all $k_0 \geq k_{\mathcal Q}$.
\end{lemma}
%%%%%%%%%%%%%%%%%%%%%%%%%%%%%%%%%%%%%%%%%%%%%%%%%%%%%%%
Our second main result of this paper is  to
characterize $q_{k_0}(x)$ in terms of the
 frequency band $(0, k_0)$, and to show that  the recovery of  $q_{k_0}(x)$ is not sensitive
 to errors in the measurements. \\

We are now ready to give the proof of Theorem~\ref{estimateobservable}.\\[1ex]

\begin{proof}
Let $p_{k_0, \pm} (x,k) $ and $\tilde p_{k_0, \pm} (x,k) $ be the impedance
functions solutions to the
system~\eqref{truncated}-\eqref{tracebis}-\eqref{boundarytruncated}
related respectively to the observable mediums $q_{k_0}$ and $\tilde q_{k_0}$.
To simplify the notation  we introduce the impedance perturbations
 $u_\pm(x,k) =  p_{k_0, \pm} (x,k) - \tilde p_{k_0, \pm} (x,k) $  due to the
measurements difference on the boundary $\epsilon(k) = d_+(k) -\tilde d_+(k)$.\\

Then $u_\pm(x,k), q_{k_0} $ and $\tilde q_{k_0}$ verify
\bean  \label{uplus}
 u_{ +}^\prime + ik (p_{k_0,+} +\tilde p_{k_0,+} ) - ik ( q_{k_0}- \tilde q_{k_0}) &=&0,\\
\label{uminus}
u_{ -}^\prime - ik (p_{k_0,-} +\tilde p_{k_0,-} )+ ik ( q_{k_0}- \tilde q_{k_0})&=&0,\\
\label{Qequation}
\left(\log \left| \frac{1+q_{k_0}}{ 1+\tilde q_{k_0}}\right|\right)^\prime-
\frac{2}{\pi}\int_{-k_0}^{k_0}( u_{+}(x,k)-u_{-}(x,k))dk &=&0,
\eean
 subject to the boundary conditions
\bean
\ba{lllccccccc} \label{bound}
u_+(0,k) &=& \epsilon(k);\; \; &   u_{ -}(0,k) &=& 0;\;\;  & q_{k_0}(0) = \tilde q_{k_0} =  0.
\ea
\eean
for all $x\in (0,1)$, $k\in \mathbb C_+$. \\

Integrating the  equation~\eqref{Qequation}  over $(0,x)$, we obtain
\bean \label{newQequation}
\log \left| \frac{1+q_{k_0}}{ 1+\tilde q_{k_0}}\right| &=&
\frac{2}{\pi}\int_0^x \int_{-k_0}^{k_0}( u_{+}(t,k)-u_{-}(t,k))dk dt.
\eean
Solving the equations~\eqref{uplus} and~\eqref{uminus}  gives
\bea
u_-(x,k) \hspace{-0.3cm}&=&\hspace{-0.3cm}
-ik \int_{0}^x \hat q(t)  e^{ik \int_t^x(p_{k_0,-}(\tau,k) +\tilde p_{k_0,-}(\tau,k)) d\tau } dt\\
u_+(x,k) \hspace{-0.3cm}&=&\hspace{-0.3cm} \epsilon(k)
e^{-ik \int_0^x(p_{k_0,+}(t,k) + \tilde p_{k_0,+}(t,k)) dt }
\hspace{-0.1cm}+\hspace{-0.1cm}
ik \int_{0}^x \hat q(t)  e^{-ik \int_t^x(p_{k_0,+}(\tau,k) +\tilde p_{k_0,+}(\tau,k)) d\tau } dt,
\eea
where $\hat q(t) = q_{k_0}(t) -\tilde q_{k_0}(t) $.\\

Substituting the new expressions of $u_\pm(x,k)$ into the
equality~\eqref{newQequation}, we find

\bean \label{newnewQequation}
\log \left| \frac{1+q_{k_0}}{ 1+\tilde q_{k_0}}\right| =\\ \nonumber
 \frac{2}{\pi}\int_{-k_0}^{k_0} \epsilon(k)
e^{-ik \int_0^x(p_{k_0,+}(t,k) +\tilde p_{k_0,+}(t,k)) dt} dk +
\frac{2i}{\pi} \int_0^x \int_0^r \hat q(t) K(r,t,k_0)dt dr,
\eean
where
\bea
K(r,t,k_0) = \int_{-k_0}^{k_0} k \left( e^{-ik \int_t^r(p_{k_0,+}(\tau,k)
+\tilde p_{k_0,+}(\tau,k)) d\tau }
+e^{ik \int_t^r(p_{k_0,-}(\tau,k) +\tilde p_{k_0,-}(\tau,k)) d\tau }\right) dk,
\eea
for $r,t \in (0,1)$.
%%%%%%%%%%%%%%%%%%%%%%%%%%%%%%%%%%%
\begin{lemma}  \label{estimategamma} Under the same conditions as
 in Theorem~\ref{estimateobservable}, there exist
constants $c_{\mathcal Q}>0$
and $k_{\mathcal Q}>0$
such  that
\bea
|K(r,t,k_0)| \leq c_{\mathcal Q},
\eea
for all $r, t \in (0,1)$ and $k_0 \geq k_{\mathcal Q}$.
\end{lemma}
%%%%%%%%%%%%%%%%%%%%%%%%%%%%%%%%
\proof(Lemma \ref{estimategamma})
First we remark from the uniqueness
of solution to the system~\eqref{truncated}-\eqref{tracebis}-\eqref{boundarytruncated}
that  $p_{k_0,+}$ like the impedance function $p_+(x,k)$,
satisfies
\bea
\overline{p_{k_0,+}(x,k)} &=&p_{k_0,+}(x,-k),
\eea
for all $x\in (0,1)$. Then, by a change of variables ($k \rightarrow -k$), we obtain
\bea
\int_{-k_0}^{k_0} k  e^{-ik \int_t^r(p_{k_0,+}(\tau,k)
+\tilde p_{k_0,+}(\tau,k)) d\tau } dk &=& -\int_{-k_0}^{k_0} k
  e^{ik \int_t^r(\overline{p_{k_0,+}(\tau,k)}
+\overline{\tilde p_{k_0,+}(\tau,k))} d\tau } dk
\eea
Hence, $K$ can be rewritten as
\bea
K(r,t,k_0) = -\int_{-k_0}^{k_0} k \left( e^{ik \int_t^r(\overline{p_{k_0,+}(\tau,k)}
+\overline{\tilde p_{k_0,+}(\tau,k))} d\tau }
-e^{ik \int_t^r(p_{k_0,-}(\tau,k) +\tilde p_{k_0,-}(\tau,k)) d\tau }\right) dk,
\eea
Now, let  $\widetilde K$ be defined as follows
\bea
K(r,t,k_0) = -\int_{-k_0}^{k_0} k \left( e^{ik \int_t^r( \overline{p_{+}(\tau,k)}
+\overline{\tilde p_{+}(\tau,k))} d\tau }
-e^{ik \int_t^r(p_{-}(\tau,k) +\tilde p_{-}(\tau,k)) d\tau }\right) dk,
\eea
According to Lemma~\ref{qpapproximation}, there
the integrand of  $K(r,t,k_0)  - \widetilde K(r,t,k_0)$
decays like $\frac{1}{k_0^{m-1}}$ uniformly with respect
to $r, t \in [0,1]$. Therefore there  exist constants
$c_{\mathcal Q} >0$ and $k_{\mathcal Q}>0$ such that
\bea
\left |K(r,t,k_0)  - \widetilde K(r,t,k_0)\right | \leq
c_{\mathcal Q},
\eea
for all $k \geq k_{\mathcal Q}$. \\

The asymptotic expansions~\eqref{asympp}
and~\eqref{asymdiffIpm} in Theorem~\ref{thmppm}
 imply that
\bea
\left|e^{ik \int_t^r \overline{p_{\pm}(\tau,k)} d\tau } \right|,
\left|e^{ik \int_t^r \overline{p_{\pm}(\tau,k)} d\tau} \right| \leq c_{\mathcal Q}
\eea
for all $t,r \in [0,1]$ and $k\in \mathbb C^+$. Furthermore
\bea
\left|e^{ik \int_t^r \overline{p_{\pm}(\tau,k)} d\tau} -
e^{ik \int_t^r \overline{p_{\pm}(\tau,k)} d\tau} \right| \leq \frac{c_{\mathcal Q}}{k_0^{m}}
\eea
all $t,r \in [0,1]$ and $k_0\geq k_{\mathcal Q}$. Combining the previous
inequalities we finally obtain that $\widetilde K (r,t,k_0)$ is uniformly
bounded over $[0,1]^2$ for all  $k_0\geq k_{\mathcal Q}$, which finishes the
proof of the lemma.\eproof
\end{proof}

Back to the equation~\eqref{newnewQequation}, by combining the
integral equation with the estimates of Lemma~\ref{estimategamma}
and the bounds over the functions $p_{k_0,+}$ and $\overline p_{k_0,+}$,
we obtain
\bean
\left| \log \left| \frac{1+q_{k_0}}{ 1+\tilde q_{k_0}}\right|\right| \leq c_{\mathcal Q}
\left( \|\epsilon(k)\|_{L^1(-k_0,k_0)} +
\int_0^x \int_0^r |\hat q(t)| dt dr \right), \nonumber\\
\leq c_{\mathcal Q}
\left( \|\epsilon(k)\|_{L^1(-k_0,k_0)} +
\int_0^x  |\hat q(t)| dt \right),
 \label{inequalityQ1}
\eean
for all $x\in (0,1)$. \\

Observing that the fact that $q_{k_0} \rightarrow q$ and $\tilde q_{k_0} \rightarrow \tilde q$
in $L^\infty(0,1)$ combined with inequalities~\eqref{refractivebounds}
imply that
 the  functions $ 1+q_{k_0}$ and $ 1+\tilde q_{k_0}$ are lower
ad upper bounded for
large $k_0$, that is, there exist a constant $k_{\mathcal Q}>0$ such that
\bea
\frac{n_{0}}{2}\leq 1+q_{k_0}(x), \; 1+\tilde q_{k_0} \leq 2n_0
\eea
for all $x \in [0,1]$ and $k_0 \geq k_{\mathcal Q}$.
Therefore
\bea
\hat q(x) \leq \frac{1}{2n_0} \left| \log \left| \frac{1+q_{k_0}}{ 1+\tilde q_{k_0}}\right|\right|,
\eea
for all $x\in [0,1]$. Combining the last inequality with~\eqref{inequalityQ1}
gives
\bean
|\hat q(x)|
\leq c_{\mathcal Q}
\left( \|\epsilon(k)\|_{L^1(0,k_0)} +
\int_0^x  |\hat q(t)| dt \right),
 \label{inequalityQ2}
\eean
for all $x\in (0,1)$ and $k_0 \geq k_{\mathcal Q}$. \\

Applying  Gronwall's inequality (Lemma~\ref{Gronwall}) on~\eqref{inequalityQ2},  with the choice
of $\rho_{\mathcal Q} = c_{\mathcal Q} +c_{\mathcal Q}^2 e^{c_{\mathcal Q}}$, we find
\bea
|\hat q(x)| \leq \rho_{\mathcal Q} \|\epsilon(k)\|_{L^1(0,k_0)}  \label{inequalityQ3}
\eea
for all $x\in \mathbb R$ and $k_0 \geq k_{\mathcal Q}$, which finishes
the proof of the Theorem ~\ref{estimateobservable}.
\eproof

\begin{remark}
%{\bf Remark 3.1.}
The estimate of Theorem \ref{estimateobservable} provides a basis for excellent numerical results
to reconstruct the observable part of the medium. In addition, it is an integral part of the proof of Theorem \ref{mainhigh}.
\end{remark}

Now, we go back to the proof of the main theorems. Lemma~\ref{qpapproximation}
implies that  if $k_0$ is
large enough we have the existence of $q_{k_0}$ and  $\tilde q_{k_0}$. By splitting
the difference $q-\tilde q$ into three parts we have
\bea
\| q-\tilde q\|_{L^\infty(0,1)} \leq \| q-q_{k_0}\|_{L^\infty(0,1)} +
\|q_{k_0}-\tilde q_{k_0}\|_{L^\infty(0,1)} +\| \tilde q- \tilde q_{k_0}\|_{L^\infty(0,1)}.
\eea
Using now the results of  Lemma~\ref{qpapproximation}
and Theorem~\ref{estimateobservable}
to estimate each part of the right hand side  we finish the
proof of Theorem~\ref{mainhigh}.
\eproof

%%%%%%%%%%%%%%%%%%%%%%%%%%%%%%%
\begin{theorem}\label{mainhigh2}
Assume that $q, \tilde q$ be two medium
functions in  $\mathcal Q$.  Let $d_+(k)$ and $\tilde d_+(k)$
be  the boundary measurements associated respectively  to $q$
and $\tilde q$ as defined in~\eqref{dplusminus}.
Then, there exist constants
 $c_{\mathcal Q}>0$  and $k_{\mathcal Q}$
such that
\bean \label{mainhighineq}
\left\| q - \tilde q \right\|_{L^\infty(\mathbb R)} \leq c_{{\mathcal Q}}
\left( \|d_\pm-\tilde d_\pm\|_{L^1(0, k_0)} + \frac{1}{k_0^{m}}\right),
\eean
for all $k_0 \geq  k_{{\mathcal Q}} $.
\end{theorem}
Obviously this result  implies the uniqueness of the multi-frequency inverse medium,
   and  a  conditional Lipschitz stability estimate
  when the band of frequency is large enough.\\
  %%%%%%%%%%%%%%%%%%%%%%%%%%%
\begin{corollary}\label{stabhigh}  Assume that $q, \tilde q$ be two medium
functions in  $\mathcal Q$.  Let $d_+(k)$ and $\tilde d_+(k)$
be  the boundary measurements associated respectively  to $q$
and $\tilde q$ as defined in~\eqref{dplusminus}, satisfying
$\|d_\pm -\tilde d_\pm \|_{L^\infty(0,+\infty)}  <1$. Then, there exists a constant
 $c_{\mathcal Q}>0$
such that the following  Lipschitz stability
\bea
\left\| q - \tilde q \right\|_{L^\infty(\mathbb R)} \leq c_{\mathcal Q}
 \|d_+(k) -\tilde d_+(k)\|_{L^\infty(0, +\infty)}^{\frac{m}{m+1}}.
\eea
holds.
\end{corollary}
\proof

Under the same assumptions of Theorem~\ref{mainhigh2}, we have
\bean \label{aall}
\left\| q - \tilde q \right\|_{L^\infty(\mathbb R)} \leq c_{{\mathcal Q}}
\left( k_0\|d_+ -\tilde d_+\|_{L^\infty(0, k_0)} + \frac{1}{k_0^{m}}\right),
\eean
for all $k_0= s k_{\mathcal Q}$ with $s>1$.
By taking $s= \|d_+ -\tilde d_+\|_{L^\infty(0, k_0)}^{-\frac{1}{m+1}}$,
we
get the wanted estimate.

\eproof

%%%%%%%%%%%%%%%%%%%%%%%%%%%%%%%%%%

%%%%%%%%%%%%%%%%%%%%%%%%%%%%%%%%%%
\begin{remark} The estimate~\eqref{mainhighineq}
 has two parts: the first is Lipschitz in terms
of the errors in measurements,  and the second
decays as the size of the frequency interval
takes larger values.  Clearly, this shows that as
the frequency increases a conditional H\"older
stability in $L^\infty $ norm can be reached as illustrated
in Corollary~\ref{stabhigh}.

\end{remark}
%%%%%%%%%%%%%%%%%%%%%%%%%%%

%%%%%%%%%%%%%%%%%%%%%%%%%%%%%%%%%%%%%
%%%%%%%%%%%%%%%%%%%%%%%%%%%%%%%%%%%
\section{Proof of Theorem~\ref{mainhigh}} \label{proofmai}
In this section we prove the stability estimate \eqref{mainestimate}.
We first provide the following  conditional stability estimate for the
unique continuation of
$d_\pm$  on a line.
%%%%%%%%%%%%%%%
\begin{theorem}\label{uniquecontinuation}
Let $k_0>0$, $d_\pm$ and $\tilde d_\pm$ be the impedance coefficients given in
 \eqref{impedanceboundaryconditions} for respectively $q$ and
 $\tilde q$ in $\mathcal Q$..
 Then
 the following estimate hold
 \bean
 |d_\pm -\tilde d_\pm |(k) \leq 2d_{\mathcal Q}
 \|d_\pm -\tilde d_\pm \|_{L^\infty(0,k_0)}^{w_0(k, k_0)},
 \eean
 for all $k\geq k_0$, where $d_{\mathcal q}$ is the constant appearing in Proposition~\ref{upperd}.
\end{theorem}
%%%%%%%%%%%%%%%%%
\proof

We deduce from Proposition~\ref{upperd} that
\bean \label{ppp2}
|d_-(k) -\tilde d_-(k) | \leq 2 d_{\mathcal Q},
\eean
for all $k \in S_{\mathcal Q}$.  \\

Without loss of generality we can assume that $h_{\mathcal Q} = \frac{\pi}{2n_{\mathcal Q}}, $ where
$n_{\mathcal Q}\in \mathbb N^*$. Let  $S_{h_{\mathcal Q}} = \{k\in \mathbb C; \re(k)>0, \, |\im(k)| <  h_{\mathcal Q}\}$,
be half a strip, and
et $w_0(k; k_0)$ be  the harmonic measure of the complex open domain
$S_{h_{\mathcal Q}} \setminus [0, k_0]\times \{0\}$. It is the unique solution  to the system:
\bea
\Delta w(k; k_0)&=&0 \quad k\in S_{h_{\mathcal Q}} \setminus [0, k_0]\times\{0\},\\
w(k; k_0) &=& 0 \quad k\in \partial S_{h_{\mathcal Q}},\\
w(k; k_0) &=& 1 \quad k\in (0, k_0]\times\{0\}.
\eea
The holomorphic unique continuation  of  the functions $d_\pm -\tilde d_\pm$ using the Two constants
Theorem \cite{Is, Ne}, gives
 \bea
  \|d_\pm -\tilde d_\pm \|_{L^\infty(0,k)} \leq (2d_{\mathcal Q})^{1-w_0(k, k_0)}
 \|d_\pm -\tilde d_\pm \|_{L^\infty(0,k_0)}^{w_0(k, k_0)}, \qquad \forall k\geq k_0.
 \eea

Finally, the bounds  satisfied by $w(k; k_0) $ are obtained from
 Lemma~\ref{harmonicmeasure}.

  \eproof

 We deduce again from  from Proposition~\ref{upperd} the existence of $k^\star \in \mathbb R_+$ satisfying
 \bea
\|d_\pm -\tilde d_\pm \|_{L^\infty(0,+\infty)}  = |d_-(k^\star) -\tilde d_-(k^\star)|.
\eea

We then deduce from Theorem \ref{uniquecontinuation}  the following estimate
\bea
\|d_\pm -\tilde d_\pm \|_{L^\infty(0,+\infty)} =|d_-(k^\star) -\tilde d_-(k^\star)| \leq 2d_{\mathcal Q}
 \|d_\pm -\tilde d_\pm \|_{L^\infty(0,k_0)}^{w_0(k^\star, k_0)}.
\eea
Considering  the global stability estimate in Corollary \ref{stabhigh}, we obtain
\bea
\left\| q - \tilde q \right\|_{L^\infty(\mathbb R)} \leq c_{\mathcal Q}
 \|d_+ -\tilde d_+\|_{L^\infty(0, +\infty)}^{\frac{m}{m+1}} \leq 2c_{\mathcal Q} d_{\mathcal Q}
 \|d_\pm -\tilde d_\pm \|_{L^\infty(0,k_0)}^{ \frac{m}{m+1}w_0(k^\star, k_0)},
\eea
which finishes the proof of the theorem.

 %%%%%%%%%%%%%%%%%%%%%%%%%%%%%%%%%%%%%%%%%%%%
 \section{Proof of Theorem~\ref{mainhigh3}} \label{proofmai2}
 In this section we prove the stability estimates \eqref{mainestimate2}-\eqref{mainestimate3}. We start
 by deriving  a lower bound to the harmonic measure $w_0$ on $\mathbb R_+$.

 \begin{proposition} \label{lowerestimatew0}
 The harmonic measure $w_0(k, k_0)$ satisfies
 \bea
 w_0(k, k_0) \geq \frac{6}{\pi}\eta(k_0)
 e^{- n_{\mathcal Q} k},
 \eea
 \end{proposition}
 \proof

It is known in the literature that  the following inequality \cite{Sha}
 \bea
 \arctan(x) \geq 3\hat \eta(x),
 \eea
 holds for all $x>0$, where
 \bea
 \hat \eta(x) = \frac{x}{1+2\sqrt{1+x^2}}.
 \eea

 Hence
 \bea
 \frac{2}{\pi} \arctan(\frac{ (e^{k_0}-1)^{{n_{\mathcal Q}}}}{\sqrt{(e^k-1)^{{2n_{\mathcal Q}}}
 - (e^{k_0}-1)^{{2n_{\mathcal Q}}}}}) \geq
 \frac{2}{\pi} \arctan\left((e^{k_0}-1)^{{n_{\mathcal Q}}} e^{-n_{\mathcal Q} k}\right)
 \geq \frac{6}{\pi} \eta(k_0)
 e^{- n_{\mathcal Q} k},
 \eea
 where $\eta(k_0) = \hat \eta ((e^{k_0}-1)^{{n_{\mathcal Q}}})$.
 \eproof

 We deduce  from Proposition~\ref{thmppm} that
\bean \label{nmb1}
|d_+(k) -\tilde d_+(k) | \leq  |d_+(k) - 1 |  + |\tilde d_+(k) -1 |  \leq   \frac{c_{\mathcal Q}} {k^m},
\eean
for all $k \in \mathbb R_+^*$, with $ c_{\mathcal Q} \geq 2d_{\mathcal Q}$. \\

Theorem \ref{uniquecontinuation}   and the last inequality lead to
 \bea
 |d_+(k) -\tilde d_+(k) |  \leq
 \min\{ 2d_{\mathcal Q}\varepsilon^{w_0(k, k_0)};   \frac{ c_{\mathcal Q}} {k^m}\},
 \eea
for all $k \in \mathbb R_+^*$.\\

Now we consider the two following cases. \\

 Case 1:  assume that $\frac{c_{\mathcal Q, 1}} {k_0^m} \leq  \varepsilon$ holds. \\

 Hence
 $\|d_\pm -\tilde d_\pm \|_{L^\infty(0,+\infty)}  \leq   \varepsilon$  is satisfied, and
we immediately get the first stability estimate  \eqref{mainestimate2}.\\

Case 2: assume that $\frac{c_{\mathcal Q}} {k_0^m}> \varepsilon$ holds. Due to the monotonicity of the functions
$ w_0(k_1, k_0)$ and $\frac{1}{k^m}$,  there exists a unique
$k_1\in (k_0, +\infty)$ satisfying
\bean \label{eeq}
\frac{c_{\mathcal Q}} {k_1^m} =  2d_{\mathcal Q}\varepsilon^{w_0(k_1, k_0)},
\eean

and

\bean \label{finalineqq}
\|d_+ -\tilde d_+\|_{L^\infty(0, +\infty)} \leq \frac{c_{\mathcal Q}} {k_1^m}.
\eean
Since $0<\varepsilon <1$, and $ c_{\mathcal Q} \geq 2d_{\mathcal Q}$, we  have
$k_1>1$.\\

On the other hand combining \eqref{eeq}, and Proposition \ref{lowerestimatew0}, gives
\bea
\frac{c_{\mathcal Q}} {k_1^m}  \leq \varepsilon^{\frac{6}{\pi}\eta(k_0)
 e^{- n_{\mathcal Q} k}},
\eea
which in turn leads to

\bea
  e^{n_{\mathcal Q} k_1}\left(\ln(2d_{\mathcal Q})- \ln(c_{\mathcal Q}) +m\ln(k_1) \right)
   \geq \frac{6}{\pi}\eta(k_0)  | \ln(\varepsilon)|.
\eea
Since $ c_{\mathcal Q} \geq 2d_{\mathcal Q}$, and
$k_1>1$, we deduce from the last inequality the existence of $c_{\mathcal Q}>0$ such that

\bea
  e^{c_{\mathcal Q}  k_1} \geq \eta(k_0)  | \ln(\varepsilon)|,
\eea
holds. Hence
 \bea
 k_1c_{\mathcal Q}   \geq  \ln\left(\eta(k_0)  | \ln(\varepsilon)|\right).
  \eea
 Combining now the last inequality and  estimate \eqref{finalineqq}, we find
 \bea
\|d_+ -\tilde d_+\|_{L^\infty(0, +\infty)} \leq \frac{c_{\mathcal Q}} {\left(\ln \left(\eta(k_0)  | \ln(\varepsilon)|\right) \right)^m }.
 \eea
By  Corollary~\ref{stabhigh}, and the last inequality, we
obtain the desired stability estimate \eqref{mainestimate3}, with $k_{\mathcal Q}= c_{\mathcal Q}^{\frac{m}{m+1}}$.

%%%%%%%%%%%%%%%%%%%%%%%%%%%%%%%%%%%%%
%%%%%%%%%%%%%%%%%%%%%%%%%%%%%%%%%%%
\section{Appendix}
We first  recall the Gornwall's inequality.
%%%%%%%%%%%%%%%%%%%%%%%%%%%%
\begin{lemma} \label{Gronwall} Assume that $u, v$ and $w: [0,1] \rightarrow \mathbb R_+$
 are  continuous functions
satisfying the inequality
\bea
u(x) \leq v(x) +\int_0^x u(t) w(t) dt,
\eea
for all $x\in [0,1]$. Then
\bea
u(x) \leq v(x) +\int_0^x v(t) w(t) e^{\int_t^xw(\tau) d\tau} dt.
\eea
\end{lemma}
%%%%%%%%%%%%%%%%%%

We  next give upper and lower estimates of a harmonic measure in a  complex strip containing a slit.
%%%%%%%%%%%%%%%%%%%%%%%%%%%%
\begin{lemma} \label{harmonicmeasure} Fix $n^\star \in \mathbb N^*$, and
let $h^\star= \frac{\pi}{2n^\star} $,  $k_0>0$ be two  fixed real constants,
$S_{h^\star} = \{k\in \mathbb C; \re(k)>0, \, |\im(k)| <  h^\star\}$ be half a strip.
Denote $w_0(k, k_0)$    the harmonic measure
of $S_{h^\star}\setminus (0, k_0]\times\{0\}$, Then

\bea
\frac{2}{\pi} \arctan(\frac{ (e^{k_0}-1)^{{n^\star}}}{\sqrt{(e^k-1)^{{2n^\star}} - (e^{k_0}-1)^{{2n^\star}}}})
\leq w_0(k, k_0) \leq  \frac{2}{\pi} \arctan\left(\inf\{\frac{k_0}{\sqrt{k^2- k_0^2}}, \frac{ e^{k_0{n^\star}}}{\sqrt{e^{2k{n^\star}}
- e^{2k_0{n^\star}}} }\}\right),
\eea
for all $k \geq k_0$.
\end{lemma}
%%%%%%%%%%%%%%%%%%%%%%%%%%%

\proof
For $n \in \mathbb N^*$,  denote by  $w_{n}(k, k_0)$
the harmonic measure of  $[0, k_0]\times\{0\}$ in
 the sector $\mathbb S_{\frac{\pi}{2n}} = \{k\in \mathbb C;  |\arg(k)| <
 \frac{\pi}{2n})  \} $. \\

 Let  $\Xi_n(k, k_0) = \sqrt{k^{2n}- k_0^{2n}}$ be  the conformal mapping  of the domain
$ \mathbb S_{\frac{\pi}{2n}}\setminus [0, k_0]\times\{0\}$ onto the right half-plane $ \mathbb S_{\frac{\pi}{2}}$. Here
 $\sqrt{k}$ is the principal branch of square root function on $\mathbb C \setminus (-\infty, 0)$ satisfying   $\sqrt{1} = 1$.
 The  parts of the   boundary  $[0, k_0]\times\{0\}|_\pm$ are then mapped onto $[-i k_0^n, i k_0^n ]$. \\

 Now
define  $w^\star(z, k_0^n)$ to  be the
harmonic measure of the right half-plane  $ \mathbb S_{\frac{\pi}{2}} \setminus [-i k_0^n, i k_0^n ]$.
The explicit expression of  $w^*$ is well known \cite{Ga}
\bea
w^*(z, k_0^n)= \frac{2}{\pi} \arctan(\frac{k_0^n}{z}), \textrm{ for }  z\in (0, +\infty).
\eea
Since $w_n(k, k_0) = w^*(\Xi_n(k, k_0), k_0^n)$ for $k\in  \mathbb S_{\frac{\pi}{2n}}\setminus [0, k_0]\times\{0\}$, we also
obtain
\bean
w_n(k, k_0)= \frac{2}{\pi} \arctan(\frac{k_0^n}{\sqrt{k^{2n}- k_0^{2n}}}), \textrm{ for }  k \in (k_0, +\infty).
\eean

 Let  $ \Xi_{-1}(k) =  e^k, $ be  the conformal mapping  of the domain
$S_{\frac{\pi}{{2n^\star}}}\setminus [0, k_0]\times\{0\}$ onto the domain $ \mathbb S_{\frac{\pi}{2n^\star}}   \setminus  \overline{B_1(0)}
\cup\left(  [1, e^{k_0}]\times \{0\}\right) $.\\

Since  $w_0(k, k_0) \leq  w_{{n^\star}}(\Xi_{-1}(k), \Xi_{-1}(k_0))$  on $ \partial \left( \mathbb S_{\frac{\pi}{2n^\star}}   \setminus  \overline{B_1(0)}
\cup\left(  [1, e^{k_0}]\times \{0\}\right) \right) $, we deduce from   the maximum principle
\bean \label{rightine1}
w_0(k, k_0) \leq  \frac{2}{\pi} \arctan(\frac{ e^{k_0{n^\star}}}{\sqrt{e^{2k{n^\star}} - e^{2k_0{n^\star}}} }),
\eean
 for all
$k \geq k_0$.
By construction we have $ S_{\frac{\pi}{2n^\star} }  \subset \mathbb S_{\frac{\pi}{2}}$, and consequently $ 0=w_0(k, k_0)
 \leq w_{2}(k, k_0)$
on $\{|\im(k)| = \frac{\pi}{2n^\star}  \}$. Then again by  the maximum   principle   we obtain
\bean \label{rightine2}
w_{0}(k, k_0) \leq w_{1}(k, k_0)=  \frac{2}{\pi} \arctan(\frac{k_0}{\sqrt{k^2- k_0^2}}),
\eean
 for all
$k \geq k_0$.\\

Combining inequalities \eqref{rightine1} and \eqref{rightine2}, we finally find
\bean \label{rightine}
w_{0}(k, k_0) \leq \frac{2}{\pi} \arctan\left(\inf\{\frac{k_0}{\sqrt{k^2- k_0^2}}, \frac{ e^{k_0{n^\star}}}{\sqrt{e^{2k{n^\star}}
- e^{2k_0{n^\star}}} }\}\right),
\eean
which gives the right-hand side inequality.\\

 Let  $ \Xi_{-2}(k) =  e^k - 1,$ be  the conformal mapping  of the domain
$ \mathbb S_{\frac{\pi}{n}}\setminus [0, k_0]\times\{0\}$ onto the domain $ D_{{n^\star}}   \setminus [0, e^{k_0}-1]\times \{0\} $,
where $D_{{n^\star}}= \{z\in \mathbb C; z+1 \in \mathbb S_{\frac{\pi}{{2n^\star}}}, \re(z)+1>0\}$.
Then $w_0(\Xi^{-1}_{-2}(k), k_0)$ is the harmonic measure of $ [0, e^{k_0}-1]\times \{0\} $ in the domain
$D_{{n^\star}}$. Now since $ [0, e^{k_0}-1]\times \{0\}  \subset \mathbb S_{\frac{\pi}{{2n^\star}}} \subset D_{{n^\star}}$, we have
$ 0=w_{{n^\star}}(k, k_0) \leq w_0(\Xi^{-1}_{-2}(k), k_0)$ on $\partial \mathbb S_{\frac{\pi}{{2n^\star}}} $. \\

The maximum principle implies that $w_{{n^\star}}(k, e^{k_0}-1) \leq w_0(\Xi^{-1}_{-2}(k), k_0)$ holds on $\mathbb S_{\frac{\pi}{{2n^\star}}} $,
and particularly, we  have
\bea
\frac{2}{\pi} \arctan(\frac{(e^{k_0}-1)^{{n^\star}}}{\sqrt{ k^{{2n^\star}}- (e^{k_0}-1)^{{2n^\star}} }})
\leq w_0(\Xi^{-1}_{-2}(k), k_0), \textrm{ for all } k \in (e^{k_0}-1, +\infty),
\eea
or equivalently
\bean \label{eqfinal}
\frac{2}{\pi} \arctan(\frac{(e^{k_0}-1)^{{n^\star}}}{\sqrt{ (e^{k}-1)^{{2n^\star}}- (e^{k_0}-1)^{{2n^\star} } }})
\leq w_0(k, k_0), \textrm{ for all } k \in (k_0, +\infty),
\eean
which provides the  desired  left-hand inequality.

\eproof

%%%%%%%%%%%%%%%%%%%%%%%%
\section*{Acknowledgements}
The work of GB was supported in part by a NSFC Innovative Group Fund
(No.11621101). The work of
FT was supported by the grant ANR-17-CE40-0029 of the French National Research Agency ANR (project MultiOnde).

%%%%%%%%%%%%%%%%%%%%%%%%%%%%%%%%%%%%%%%%%%%?%%%%%%%%%%%%%%%%%%%%%%%%%%%%%%%%%%%%%%%%%

%\bibliographystyle{abbrv} % apa abbr
%\bibliography{\latexpath/main}
\end{document}